\documentclass[11pt]{article}
\DeclareSymbolFont{AMSb}{U}{msb}{m}{n}
\DeclareMathSymbol{\N}{\mathbin}{AMSb}{"4E}
\DeclareMathSymbol{\Z}{\mathbin}{AMSb}{"5A}
\DeclareMathSymbol{\R}{\mathbin}{AMSb}{"52}
\DeclareMathSymbol{\Q}{\mathbin}{AMSb}{"51}
\DeclareMathSymbol{\I}{\mathbin}{AMSb}{"49}
\DeclareMathSymbol{\C}{\mathbin}{AMSb}{"43}
\def\qed{\quad \vrule height7.5pt width4.17pt depth0pt}

\begin{document}

      	\newtheorem{definition}{Definition}
	\newtheorem{theorem}[definition]{Theorem}
	\newtheorem{lemma}[definition]{Lemma}
	\newtheorem{corollary}[definition]{Corollary}
	\newtheorem{proposition}[definition]{Proposition}

\begin{center}
\textbf{Dehn Filling and Asymptotically Hyperbolic Einstein
Manifolds}

\medskip
Gordon Craig\\
Bishop's University\\
Lennoxville, Quebec, Canada\\
gcraig@ubishops.ca\\
\end{center}
\medskip

\begin{abstract}
In this article, we extend Anderson's higher-dimensional Dehn
filling construction to a large class of infinite-volume
hyperbolic manifolds. This gives an infinite family of
topologically distinct asymptotically hyperbolic Einstein
manifolds with the same conformal infinity. The construction involves
finding
a sequence of approximate solutions to the Einstein equations and then
perturbing them to exact ones.
\end{abstract}

	\section{Introduction}

In this article we will describe a construction to generate infinitely 
many non-homotopic asymptotically hyperbolic Einstein metrics with a 
shared conformal infinity. These are obtained by ``capping off'', or filling in the cusps 
of certain complete hyperbolic manifolds. The construction is based on 
Anderson's Dehn filling result for finite-volume hyperbolic manifolds(\cite{1}), and 
we will refer readers to that paper where it is appropriate. 

We will begin by providing background to this result and discussing its significance. After a brief sketch of the proof, the rest of this paper is devoted to the details.
	
Roughly speaking, asymptotically hyperbolic(AH) manifolds are complete 
Riemannian manifolds with an ideal boundary at infinity. The canonical 
example is the Poincar\'e sphere with its sphere at infinity. For this to 
work, the underlying differential manifold must be a smooth, compact 
manifold with boundary.
	
	\begin{definition}
	Let $\bar{M}$ be a compact manifold with nonempty boundary. A 
smooth 
function $\rho$ is said to be a defining function for $\partial M$ iff
\begin{equation}	
\rho : \bar{M} \longrightarrow [0,\infty)
\end{equation}
satisfies $\rho(p)=0$ iff $p \in \partial M$ and $d\rho \neq 0$ on $\partial M$.
\end{definition}

Then we have:

\begin{definition}
A complete metric $g$ on $M = \mbox{int}(\bar{M})$ is said to be conformally compact iff there exists a defining function 
$\rho$ for $\partial M$ such that 
$\bar{g}= \rho^{2}g$ extends to a metric on $\bar{M}$.
\end{definition}

As the name suggests, AH manifolds have curvature decaying to $-1$ at 
infinity. This is almost equivalent to conformal compactness, but we need 
assume a small amount of regularity on the boundary. In what follows, 
quantities with a bar over them will be measured with
respect to the compactified metric $\bar{g}$, and ones without bars will be measured with respect to the metric $g$ on $M=int(\bar{M})$. 

	\begin{definition}
	Consider a conformally compact metric $g$ on $M$. If there exists a defining function $\rho$ such
that $|\bar{\nabla} \rho| \equiv 1$ on $\partial M$, then we say 
that $(M,g)$ is asymptotically hyperbolic.
	\end{definition}

The reason for this terminology is that in this case the sectional 
curvatures of $(M,g)$ tend uniformly to -1(c.f. \cite{3}.)

	Asymptotically hyperbolic manifolds are a natural class of 
non-compact manifolds to work with because they have a nice structure at 
infinity; their curvature is tending toward a constant, and via the 
compactification $\bar{g}$ they have a ``boundary 
metric'' at infinity. Since this metric is determined by the choice of 
the function $\rho$, it only makes sense to speak of a boundary conformal class. This conformal class is known as 
the conformal infinity of the complete manifold $(M, g)$. 
	
	Given an AH manifold, it is natural to want to put a canonical AH metric on it. In two or three dimensions, the natural choice is a hyperbolic 
metric. In higher dimensions, however, hyperbolic metrics generalize in two
ways: to hyperbolic metrics and to negatively curved Einstein metrics(constant negative Ricci curvature.) If $n>3$, the
curvature tensor has more components than the metric, so prescribing sectional curvature becomes more difficult. On the
other hand, the Ricci tensor has the same rank as the metric, so Einstein metrics are natural candidates to be canonical from an algebraic/PDE point of view. Accordingly, we define:

	\begin{definition}
	An asymptotically hyperbolic manifold with constant Ricci 
curvature is called an asymptotically hyperbolic Einstein(AHE) manifold.
	\end{definition}

	AHE manifolds are sometimes referred to as Poincar\'e-Einstein manifolds. There 
has been a great deal of interest in AHE metrics recently 
due to applications to physics(c.f. \cite{27}, \cite{28}, \cite{25}.) Physicists are particularly interested in the 
correspondence between AHE metrics and their conformal
infinities. From a more analytic/geometric perspective, this 
correspondence can be thought of as a geometric 
Dirichlet problem(although a priori the topology of the filling manifold could
be undetermined.) 

	Anderson has worked extensively on this correspondence in dimension 4(c.f. the surveys \cite{28}, \cite{79} and for more details \cite{2}, \cite{3}). 
In this case, there is generally not a bijective correspondence between 
AHE metrics and their conformal infinities, but under certain geometric conditions on the boundary, there are only finitely
many AHE manifolds bounded by a given conformal class. This is not always the case, however; in \cite{2}, Anderson constructs infinitely many 
AHE metrics bounded by the same conformal class. He then shows in \cite{2} that in dimension 4, any such collection has a 
limit point which is an AHE manifold with cusps(i.e. an Einstein metric 
whose ends are either conformally compact or of finite volume.) Furthermore, this limit point has the same 
conformal infinity as all the elements in the set. 

	Under some fairly natural conditions on these limit manifolds, it is possible to show(still in 4 dimensions) that these limit manifolds are 
actually hyperbolic(see \cite{2}.) This suggests a natural question: given a hyperbolic manifold whose ends are either conformally compact or cusps, is there a sequence of AHE metrics with the same same conformal infinity converging toward it? This is indeed the case in three dimensions(where AHE metrics are 
hyperbolic,) although the methods used in this case come from hyperbolic geometry, and cannot be applied to Einstein manifolds. Anderson has recently developed a cusp closing technique for finite-volume 
hyperbolic manifolds, generalizing Thurston's Dehn filling result to higher dimensions(c.f. \cite{1}.) His construction leads to infinite families of 
topologically distinct compact Einstein 
manifolds.

	Our main result applies Anderson's cusp closing construction to generate a host of AHE metrics with the same conformal infinity:
	
\begin{theorem} \label{main}
Let $(N^{n},g)$, $(n>2)$, be a complete geometrically finite hyperbolic manifold, all of 
whose cusps have toric cross sections.  
Then it is possible to close the cusps to obtain infinitely many metrically distinct AHE manifolds, all of which
have the same conformal infinity as the original one. If the hyperbolic manifold $N^{n}$ is nonelementary, then this 
procedure gives infinitely many homotopy types. If $n>3$ these AHE metrics are non-hyperbolic. 
\end{theorem}

It is possible to obtain a large class of manifolds satisfying the hypotheses of this theorem by taking Maskit combinations of complete hyperbolic cusps along 
hyperplanes. (c.f. \cite{23} for this, and also for some background on nonelementary and geometrically finite hyperbolic manifolds.) In dimension 3, our 
construction is simply a PDE proof of Thurston's Dehn-filling theorem, and so the manifolds 
we obtain are already well-known(c.f. \cite{68}.) In higher dimensions, 
however, it gives new examples of complete negatively curved Einstein manifolds. 

From the point of view of the Dirichlet problem for AHE manifolds, we may restate our main result as follows:

\begin{proposition}
Let the conformal class [$\gamma$] be the conformal infinity of a nonelementary geometrically finite hyperbolic manifold, whose ends are all either expanding or are toric cusps, and which has at least one cusp. Then there are infinitely many topologically distinct AHE manifolds whose conformal infinity is [$\gamma$].
\end{proposition}

It is interesting to note that all of the filling manifolds which we 
construct for $[\gamma]$ have the same Euler characteristic. If a conformal class $[\gamma]$ on a three-manifold includes a metric of positive scalar curvature, then the only way that it can bound infinitely many AHE 4-manifolds is if their Euler characteristics are unbounded. (c.f. \cite{3})

	The proof of the main result follows Anderson's proof in \cite{1}: we construct a sequence of approximate 
solutions to the Einstein equations(i.e. metrics whose Ricci curvature 
is tending toward some constant) and then perturb the metrics to an exact solution. This basic gluing procedure is quite common in geometric analysis, c.f. for example \cite{13}, \cite{8}, \cite{7}, \cite{9}.

The choice of approximate solutions is the main stumbling block 
in this procedure. Although in theory it is easy to prescribe Ricci curvature, since the Ricci tensor is of the same rank as the metric, in practice it is 
extremely difficult, because we must find explicit solutions to a coupled system of nonlinear PDEs. The construction of 
the approximate solutions requires the use of a very special (explicit) family of metrics, as we shall see below.

For the perturbation argument, we will be using a functional $\Phi$, and 
metrics which satisfy $\Phi(g)=0$ will be Einstein. Then we will have a sequence of approximate solutions $g_{n}$ such 
that $\Phi(g_{n}) \longrightarrow 0$. 

It turns out that sequences of approximate solutions degenerate, so we cannot use a limiting argument to obtain our 
exact solution. On the other hand, the linearization of $\Phi$ at each metric $g_{n}$ is invertible. Thus, we could 
hope to use the inverse function theorem to invert $\Phi$ in a neighborhood of $\Phi(g_{n})$. Since $\Phi(g_{n}) 
\longrightarrow 0$, we can hope that for $n$ large enough, one of these neighborhoods will contain $0$, which will 
give 
us a metric $g$ such that $\Phi^{-1}(0)=g$. Invertibility of $\Phi$ near $\Phi(g_{n})$ is not enough to insure this 
however. We could have a situation in which the region on which $\Phi$ is invertible shrinks as $n \longrightarrow 
\infty$, so $0$ never lies in this region. 

Thus, we need to get a control on the size of the balls on which we can 
invert. If this is the case, then for $n$ large enough we can perturb 
$g_{n}$ to a metric satisfying $\Phi(g)=0$. We obtain this 
uniform control by bounding the linearization of $\Phi$.

There is no need to assume any sort of nondegeneracy hypotheses on the original hyperbolic metric. (Roughly speaking,
nondegeneracy means that a metric is a regular point for the Einstein operator. This is important if we want to apply an inverse function theorem, as we will.)
We will actually see that our approximate solutions are nondegenerate as a consequence of our main estimate (\ref{est}).

The main difference between our construction and Anderson's is that our approximate solutions are noncompact, due to the presence of the expanding ends. 
This introduces some technical difficulties in our analysis, but by 
construction we have very strong control over the expanding ends, since 
our metrics have fixed conformal infinities.

  As we mentioned above, in \cite{2}, Anderson showed that AHE manifolds with cusps which satisfy certain relatively natural
conditions are in fact hyperbolic. This leads to another interesting question: are there any non-hyperbolic AHE manifolds with cusps? Our main result allows us to construct non-hyperbolic examples by closing the cusps on a hyperbolic manifold with several
cusps, and then reopening one of them.

\begin{theorem} \label{open}
Given a nonelementary hyperbolic manifold $N^{n}$ , $n>3$, with at least one expanding end and at least two cusps with toric cross-sections, and no ends of another type, it is possible to construct infinitely many nonhomotopic nonhyperbolic Einstein manifolds which have at least one finite-volume end(i.e. a cusp) and whose infinite-volume ends are asymptotically hyperbolic. Furthermore, their expanding ends have the same conformal infinities as $N^{n}$.
\end{theorem}

Let us now set some of our notation and conventions. From now on, all manifolds will be assumed to be complete and AH, unless otherwise stated. Pointwise norms and inner products will be denoted by $|h|$ and $(f,h)$ 
respectively, while global ones will be denoted by $\|h\|$ and $\langle f,h\rangle$. $K$, $ric$, $z$ and $s$ will 
represent the sectional, Ricci, trace-free Ricci and scalar curvatures. $i(M)$ will denote the injectivity radius of $M$. $n$ will be reserved for the dimension the 
manifold $M$, and will always be strictly greater than 2. The curvature operator is defined as 
\begin{equation}
R(X,Y)Z = \nabla_{Y}\nabla_{X}Z-\nabla_{X}\nabla_{Y}Z-\nabla_{[X,Y]}Z
\end{equation}
for any three $X$, $Y$, $Z \in TM$. Our Laplacians will have negative spectrum, so $\Delta_{S^{1}}= -\frac{d^{2}}{d\theta^{2}}$. This is the so-called ``Geometer's Laplacian.'' We will often drop subscripts to improve readability if this will not lead to any confusion.

This article will be organized as follows: in section 2 we construct our approximate solutions and 
discuss some of their properties, in section 3 we will define our function spaces and operator and then in section 4 we obtain a 
uniform control 
over the operator $D\Phi$ on all the approximate solutions. Finally, in 
section 5, wrap things up by perturbing our approximate solutions to 
exact ones and discuss some points of interest, including the proof of Theorem \ref{open}.

\emph{Acknowledgements:} This work was part of the author's doctoral thesis at Stony Brook University, written under 
the supervision of Professor Michael Anderson. I would like to thank him for his guidance, advice and support. The article was 
written up during a postdoctoral fellowship at the CRM in Montreal supported in part by Professors Vestislav Apostolov 
of UQAM and Niky Kamran of McGill, and I would like to thank them for 
their support, both moral and financial. I would additionally like to express my gratitude 
to all my friends and family, and especially to Jean Dendy and Louis Garceau, for their 
encouragement. I am also very grateful to the referees for their 
comments and suggestions. Finally, I would like to dedicate this article to the memory of Ray ``Rocco" Golbert, whose support and encouragement helped make my 
time in graduate school much easier, and who will be sorely missed.

	\section{Construction of Approximate Solutions}

In this section, we construct our approximate solutions, and discuss some of their topological properties. We want to 
fill in the cusps of our hyperbolic manifold, so we will be cutting them off and filling them in with solid tori. 
Topologically, this construction is easy. Even metrically it is not too hard, assuming that we do not require anything 
of the filled manifolds. But we want our filled manifolds(our approximate solutions) to have Ricci curvature close to 
a constant. This turns out to be much more difficult. 

We will be filling each cusp separately, so we will only need to explain the procedure on one of them. By assumption, all of our 
cusps look like
\begin{equation}
g_{C}= \rho^{-2}d\rho^{2} + \rho^{2}g_{T^{n-1}}; \ \ \ \ \ \ \rho_{0} > \rho >0
\end{equation}

Note that as $\rho \rightarrow 0$, the $T^{n-1}$'s are collapsing. Without loss of generality, we can assume that 
$\rho_{0} \geq 1$ by rescaling the $\rho$ parameter. This will give us a metric of the same form, but with a rescaled 
$T^{n-1}$. Let us cut off the cusp at the torus $\rho=1$ and discard the region $0<\rho<1$. Then we are faced with the task of attaching something to 
boundary torus $T_{0}$ in such a way that the metric on the glued manifold is smooth. Note that $T_{0}$'s metric is 
the flat metric $g_{T^{n-1}}$.

For this construction to work, we will need to use a sequence of filling manifolds which are hyperbolic near their 
boundary, and whose trace-free Ricci curvature tends toward zero. We will use members of a family of AHE metrics on 
$D^{2} \times T^{n-2}$. We can obtain our filling manifolds by truncating these at some fixed distance, and then 
perturbing the metric near the boundary to make it hyperbolic. The perturbations will get smaller as we go further out, since the manifold is AH. We will start by discussing these filling manifolds.

Consider the following metric on $D^{2}\times T^{n-2}$;
\begin{equation}
g_{BH} = (V(r))^{-1}dr^{2} + V(r)d\theta^{2} + r^{2}g_{T^{n-2}}
\end{equation}
where $g_{T^{n-2}}$ is an arbitrary flat metric on the $(n-2)$-torus and we use coordinates 
$(r, \theta)$ for the disk, with $r\geq r_{+}>0$ and $\theta\in [0,\beta]$. (Note that this 
means that the locus $r=r_{+}$ is the center of the disk, so there is a coordinate 
singularity there.)

We will specify the values of the parameters $r_{+}$ and $\beta$ and the exact form of 
$V(r)$ below, but first let us 
calculate the curvatures of these metrics in terms of the function $V(r)$. We will start by setting up an orthonormal 
basis for these metrics: let $e_{1}=\sqrt{V}\partial_r$, $e_{2}=\frac{1}{\sqrt{V}}\partial_{\theta}$, and 
$e_{j}=\frac{1}{r}\partial_{\phi_{j}}$, where the $\partial_{\phi_{j}}$, $3\leq j \leq n$ are an  orthonormal basis 
for the $T^{n-2}$. A straightforward calculation shows that the $e_{i}$ 
diagonalize the curvature tensor, and that the corresponding sectional curvatures are
\begin{eqnarray} \label{curvs}
K_{12} = -\frac{V''}{2} \\
K_{1j} = K_{2j} = -\frac{V'}{2r} \ \ \ \ \ \ j>2 \\
K_{ij} = -\frac{V}{r^{2}} \ \ \ \ \ \ i,j>2 \\
\end{eqnarray} 

Now, let 
\begin{equation}
V(r)=r^{2}-2mr^{3-n}
\end{equation}

Using the same basis as above, we have
\begin{eqnarray}
K_{12} = -1+ \frac{(n-3)(n-2)m}{r^{n-1}} \\
K_{1j} = K_{2j} = -1 -\frac{(n-3)m}{r^{n-1}} \\
K_{ij} = -1 + \frac{2m}{r^{n-1}}
\end{eqnarray} 
where once again $i,j$ are assumed to be greater than 2.

Another straightforward calculation shows that this metric is Einstein with scalar curvature $-n(n-1)$ and 
asymptotically hyperbolic. (If $n=3$, it will be hyperbolic.) We 
have yet to specify the range of the $r$ parameter, but it is clear that the metric is 
well-defined for large enough $r$, so it makes sense to speak of its local and asymptotic 
properties. 

If $m=0$, we get a hyperbolic cusp metric
\begin{equation}
g_{C}= r^{-2}dr^{2} + r^{2}g_{S^{1} \times T^{n-2}}
\end{equation}
This metric will be complete if we let $r \in (0, \infty)$.

On the other hand, if $m>0$ and $n>3$, we get a nontrivial Einstein metric. These metrics 
are called $T^{n-2}$ Anti-deSitter 
Black Hole metrics. They will be complete provided we let $r\geq r_{+}= (2m)^{\frac{1}{n-1}}$ and $\theta \in 
[0,\beta_{m}]$, where 
\begin{equation}
\beta_{m}=\frac{4\pi}{(n-1)r_{+}}
\end{equation}
 (c.f. \cite{2}.) Note that the locus $\{r=r_{+}\}$ 
is a flat totally geodesic 
$T^{n-2}$. By analogy with the core geodesics in hyperbolic Dehn surgery(c.f. \cite{22},) we call this a core torus.

	Now recall that we introduced these manifolds because we want to glue them into a cusp. They have the correct 
topological and local geometric properties to work. The first choice we could make would be to cut off one of the 
black hole metrics above at some large $r$, and then perturb it to make it hyperbolic near the boundary. The problem 
is that we cannot fix the global geometry near the boundary; although we can choose the metric on the $T^{n-2}$, the 
boundary metric will necessarily be the product of this flat metric and a large $S^{1}$, since the size of the $S^{1}$ 
factor is determined by $r$.

To resolve this difficulty, we will exploit the large isometry group of these metrics to take a quotient with the 
desired boundary. Below, we shall use the term ``black hole metric'' to refer to any metric on $D^{2} \times T^{n-2}$ 
which has the same universal cover as $g_{BH}$.

\begin{proposition} \label{Quotient}
Suppose we have an  $S^{1} \times \R^{n-2}$-invariant metric on $D^{2} \times \R^{n-2}$. Let $T_{0}$ be some flat 
$(n-1)$-torus, and let $\sigma \subset T_{0}$ be a simple closed geodesic such that
\begin{equation}
L(\sigma) = L(\partial D^{2}).
\end{equation}
Then $\exists \Gamma_{0} \subset \ \mbox{Isom}(D^{2} \times \R^{n-2})$ such that
\begin{equation}
\frac{(D^{2}\times\R^{n-2})}{\Gamma_{0}} \simeq M_{fill}
\end{equation}
 is a solid torus with $\partial M_{fill} = T_{0}$.
\end{proposition}

\textbf{Proof:} 
We have 
\begin{equation}
T_{0} = \R^{n-1} / \Gamma,
\end{equation}
where $\Gamma$ is some $(n-1)$-dimensional group of translations of $\R^{n-1}$. Since $\sigma$ is closed and simple, 
the translation induced by $\sigma$ is a generator for $\Gamma$. We can find elements $\gamma_{i} \in \R^{n-1}$ such 
that the set $\{\sigma, \gamma_{1}, ..., \gamma_{n-2} \}$ forms a set of generators for $\Gamma$. Let us denote the 
subgroup of $\pi_{1}(T_{0})$ generated by $\sigma$ by $\langle \sigma \rangle$. Then since $\pi_{1}(T_{0})$ is 
Abelian,  $\langle \sigma \rangle$ is normal, which implies that the covering map
\begin{equation}
p: \R^{n-1} \longrightarrow T_{0}
\end{equation}
splits as $p = p_{1} \circ p_{2}$, where
\begin{equation}
p_{2}: \R^{n-1} \longrightarrow \R^{n-1} / \langle \sigma \rangle \cong S^{1} \times \R^{n-2}
\end{equation}
and 
\begin{equation}
p_{1}:  S^{1} \times \R^{n-2} \longrightarrow (S^{1} \times \R^{n-2})/ {\Gamma_{0}} \cong T_{0}
\end{equation}
where $\Gamma_{0} = \Gamma / \langle \sigma \rangle$

	Now, say we have an $S^{1} \times \R^{n-2}$-invariant metric on $D^{2} \times \R^{n-2}$, and that the length of 
$\partial D^{2}$ is $L(\sigma)$.

We will use the above remark to construct a quotient of this metric with boundary $T_{0}$. We can describe this 
quotient in terms of coordinates $(r, \theta, \phi_{3}, ..., \phi_{n-2})$, where the $\phi_{i}$'s are the standard 
coordinates on $\R^{n-2}$. Define an isometric action of $\Gamma_{0}$ on $D^{2} \times \R^{n-2}$ by keeping $r$ fixed 
and acting on the $S^{1} \times \R^{n-2}$ coordinates. The boundary of this quotient will certainly be $T_{0}$, and it 
is clear that there are no fixed points on $r=0$ since no element maps $\sigma$ to itself except the identity, so the 
quotient is indeed a manifold. \qed

Note that the reason that this works is that we were able to split off the $\langle \sigma \rangle$ from the rest of $\Gamma$, and then fill it in with a disk. The basic point is that when we project the other generators onto the core $\R^{n-2}$, they cannot be zero, or else they would be parallel to $\sigma$. If one thinks about the three-dimensional case, one can picture the universal cover as being a tubular neighborhood of a geodesic in hyperbolic space. Then $\sigma$ would be the boundary of a disk perpendicular to the core geodesic. One can obtain the torus $T^{2}$ by taking the quotient of the cylinder by $\langle \gamma \rangle$, where $\gamma$ is  some composition of a translation and a rotation. The only way that this action will not extend to the core geodesic is if $\gamma$ has no translation component. But this is impossible if the quotient of the boundary is to be a torus.

	Now we will get an appropriate metric on this quotient. All we need is an $S^{1} \times \R^{n-2}$-invariant 
metric on $D^{2} \times \R^{n-2}$. Since we want an Einstein metric, we will take the universal cover of the 
$T^{n-2}$-black hole metrics, slightly altered near the boundary. It turns out that the value of $m$ is irrelevant to 
the local geometry of these(c.f. \cite{2},) so we will set $m=1$.

Let \begin{equation}
R=\frac{L(\sigma)}{\beta_{1}}
\end{equation}
and define
\begin{equation}
\widetilde{g_{fill}} = V(r)^{-1} dr^{2} + V(r) d\theta^{2} + r^{2}g_{Eucl}; r \in [r_{+}, R), \theta \in [0, 
\beta_{1})
\end{equation}
where $V(r) = r^{2} - \frac{2\chi(r)}{r^{n-3}}$, $\chi(r) = 1$ for $r < R-2$, $\chi(r) = 0$ for $r > R-1$. 

Now note that for $r>R-1$, 
\begin{equation}
\widetilde{g_{fill}} = r^{-2} dr^{2} + r^{2}\left(d\theta^{2} + g_{Eucl}\right)
\end{equation}

By taking the quotient, we will get the metric

\begin{equation}
g_{fill} = r^{-2} dr^{2} + r^{2}\frac{g_{T_0}}{R^{2}}
\end{equation}
on $r>R-1$.

By the change of coordinates
\begin{equation}
\rho = \frac{r}{R}
\end{equation}
we get a metric which is identically equal near its boundary to the hyperbolic cusp metric which we are trying to fill 
in. In the event that we have $k$ toric cusps, we can cut off each one, and perform this procedure on a geodesic 
$\sigma^{i}$ in each boundary torus. We then obtain a manifold $(M_{\sigma}, g_{\sigma})$, where $\sigma$ is the 
ordered $k$-tuple of geodesics $(\sigma^{i})_{1\leq i\leq k}$.

We say that $M_{\sigma}$ is a Dehn filling of $N$, again by analogy with the 
three-dimensional case. The size of the 
Dehn filling is defined to be
\begin{equation}
|\sigma| = \min_{i} L(\sigma^{i}) = \min_{i} (R_{\sigma^{i}}\beta_{1})
\end{equation}
This is well-defined, since we fix the boundary tori at the beginning.

In three dimensions, we no longer have these black hole metrics. In fact, the only candidate for the glued-in metric 
is a quotient of hyperbolic space. But there are many hyperbolic 3-manifolds bounded by 2-tori; we can repeat the 
above Dehn filling on the tubular neighborhood of a geodesic in hyperbolic 3-space. The analysis below then goes analogously to the higher dimensional-case.

	We end this section with a few remarks on the topology of the $M_{\sigma}$'s. We start with the Gromov-Thurston 
$2\pi$-theorem. (For an elementary proof, c.f.\cite{18})
 
 	\begin{proposition} (Gromov-Thurston, \cite{19})  If $|\sigma| > 2\pi$, then $M_{\sigma}$ admits a 
nonpositively curved metric, such that the core tori are completely geodesic and whose sectional curvature is strictly 
negative off the core tori.
 	\end{proposition}

Thus, we have:

\begin{proposition} \label{Kp1}
For $|\sigma|>2\pi$, $M_{\sigma}$ is a $K(\pi, 1)$, and every noncyclic Abelian subgroup of $\pi_{1}(M_{\sigma})$ is 
carried by one of the core tori.
\end{proposition}
\textbf{Proof:} (Adapted from Theorem 6.3.9 of \cite{24}.) The first statement follows directly from the fact that 
$M_{\sigma}$ admits a metric of nonpositive sectional curvature. To prove the second one, consider the action of 
$\pi_{1}(M_{\sigma})$ by isometries on $\widetilde{M_{\sigma}}$, the universal cover of $M_{\sigma}$, where 
$\widetilde{M_{\sigma}}$ is equipped with the lifted metric. Let $\alpha\in \pi_{1}(M_{\sigma})$, and let 
$p\in\widetilde{M_{\sigma}}$ satisfy 
\begin{equation}
d(p,\alpha(p)) = \inf_{x\in \widetilde{M_{\sigma}}}d(x,\alpha(x))
\end{equation}
Such a $p$ exists because $M_{\sigma}$'s is injectivity radius is nonzero, and is large on the expanding ends of 
$M_{\sigma}$. Thus any minimizing sequence for the function $d(x,\alpha(x))$ must project down to some compact subset 
of $M_{\sigma}$.
 
Since  $\widetilde{M_{\sigma}}$ is simply connected and of nonpositive curvature, there is a unique geodesic $\gamma$ 
joining $p=\gamma(0)$ and $\alpha(p)=\gamma(t_{0})$. Furthermore, $\alpha$ must fix $\gamma$. To show this, all we 
need to prove is that $\alpha^{2}(p)$ lies on $\gamma$. If it does not, then let $x=\gamma(\frac{t_{0}}{2})$ be the 
midpoint of the segment $\overline{p\alpha(p)}$ of $\gamma$. By assumption, the three points $x$, $\alpha(x)$ and 
$\alpha(p)$ cannot all lie on the same geodesic, so we must have
\begin{eqnarray}
d(x, \alpha(x)) & < & d(x,\alpha(p)) + d(\alpha(p), \alpha(x)) \\
& = & \frac{t_{0}}{2} + d(p,x) \\
& = & t_{0} \\
& = & d(p,\alpha(p)),
\end{eqnarray}
which is impossible, since $p$ minimizes $d(x,\alpha(x))$.

If $\alpha$, $\beta \in \pi_{1}(M_{\sigma})$ commute, then $\alpha\circ\beta(p)=\beta\circ\alpha(p)$, so $\beta$ sends 
the geodesic segment $\overline{p\alpha(p)}$ to the geodesic segment $\overline{\beta(p)\beta(\alpha(p))}$. 
Furthermore, we can see that $\alpha$ must map the unique geodesic joining $p$ to $\beta(p)$ to the unique one joining 
$\alpha(p)$ and $\beta(\alpha(p))$. This gives a geodesic quadrilateral whose angles add up to $2\pi$. By the 
Topogonov theorem, this quadrilateral must lie in a flat submanifold, so we are done. \qed

We will also need the following result:

\begin{proposition} 
Let $N$ admit a nonelementary geometrically finite hyperbolic metric. Then $N$ admits only finitely many homotopy 
equivalence classes. 
\end{proposition}
\textbf{Proof:} Let 
\begin{equation}
F:N \longrightarrow N
\end{equation}
be a homotopy equivalence. Given a hyperbolic metric $g$ on $N$, we can deform $F$ so that it fixes the conformal 
infinity of $(N,g)$. Then by Sullivan rigidity(c.f. \cite{23}), $F$ can be represented by an isometry. But by 
\cite{21} generic nonelementary geometrically finite hyperbolic metrics have only finitely many isometries. \qed

We may now prove that our Dehn fillings give rise to infinitely many homotopy classes. 

\begin{proposition}
Let N admit a geometrically finite nonelementary hyperbolic metric, all of whose cusps are tori. Let $M_{\sigma}$ be 
obtained from $N$ by a Dehn filling. Then there are only finitely many other Dehn fillings which have the same 
homotopy type as $M_{\sigma}$. 
\end{proposition}
\textbf{Proof:} We will simply adapt Anderson's proof from \cite{1} to the infinite volume case. The idea of the proof 
is to show that a homotopy equivalence
\begin{equation}
F: M_{\sigma} \longrightarrow M_{\sigma^{\prime}}
\end{equation}
leads either to a nontrivial homotopy equivalence of the hyperbolic manifold $N$, or to a nontrivial isomorphism of 
the Dehn filling data. Since there are only finitely many members in either of these classes, we can conclude that 
there are only finitely many $M_{\sigma}$'s in any homotopy class.                

By Seifert-Van Kampen, we have that 
\begin{equation}
\pi_{1}(M_{\sigma}) = \frac{\pi_{1}(N)}{\langle \cup \sigma_{i} \rangle}
\end{equation}

By proposition \ref{Kp1}, we know that for $|\sigma|$ sufficiently large, the only (conjugacy classes) of noncyclic 
Abelian subgroups of $\pi_{1}(M_{\sigma})$ are carried by the core tori. Now, say that we have a homotopy equivalence
\begin{equation}
F: M_{\sigma} \longrightarrow M_{\sigma^{\prime}}
\end{equation}
Then $F_{*}$ must permute the cyclic subgroups carried by the essential tori. This in turn implies that $F$ must map 
neighborhoods of these tori onto each other. We can then use $F$ to define a map $G$ from the original hyperbolic 
manifold $N$ to itself such that $G$ fixes the conformal infinity of $N$, interchanges the cusps of $N$ and such that 
$G_{*}$ is an isomorphism of $\pi_{1}(N)$. Then again by Sullivan rigidity, $G$ is a homotopy equivalence, of which 
there are only finitely many by the previous proposition. Now the homotopy equivalence $F$ must also preserve the Dehn 
filling data, so necessarily we must have
\begin{equation}
F_{*}\langle \sigma_{i} \rangle =  \langle \sigma_{j}^{\prime} \rangle
\end{equation}
But given a cyclic group, there are only two elements which can generate it. Thus, there are only finitely many 
Dehn-filled manifolds homotopy equivalent to $M_{\sigma}$. \qed

The proposition does not hold if we drop the hypothesis that $(N,g)$ is nonelementary; one can construct infinitely 
many 
nonisometric black hole metrics on the solid torus $D^{2} \times T^{n-2}$ with the same conformal infinity 
$T^{n-1}$(c.f. \cite{2}.) These can also be thought of as Dehn fillings of the complete hyperbolic cusps
\begin{equation}
g_{C}=r^{-2}dr^{2}+r^{2}g_{T^{n-1}}
\end{equation}

\section{Analytic Preliminaries}

Let us start by establishing the function spaces we will be using. To obtain good 
analytic properties, such as elliptic 
estimates and precompactness of bounded sequences, we will work 
with H\"older spaces. The precise definition of these spaces requires some care, as we will want to be able to compare spaces defined on a large class of manifolds.

\begin{proposition} (\cite{6})
Let $Q>1$, $C, \ i_{0}>0$, $k \in \N$, $0<\alpha<1$. Then there exists $\rho_{0}>0$ such 
that if $\|\nabla^{k-1} ric\|_{L^{\infty}} \leq C$ and $i(M) \geq i_{0}$, then for any 
$x\in M$, the ball $B(x,\rho_{0})$ has
harmonic coordinates in which we have
\begin{equation}
Q^{-1}I \leq g \leq Q I 
\end{equation}
\[
\sum_{1 \leq |\beta| \leq k} \rho_{0}^{|\beta|}\sup_{y \in B(x,\rho_{0})}|\partial^{\beta}g_{ij}(y)| + \sum_{|\beta | = k}
\rho_{0}^{k+\alpha}\sup_{y_{1},y_{2} \in  B(x,\rho_{0}) }\frac{|\partial^{\beta}g_{ij}(y_{1}) - \partial^{\beta}g_{ij}(y_{2})|}{|y_{1} - y_{2}|^{\alpha}} \]
\begin{equation}
 \ \ \ \ \ \ \ \ \ \ \ \ \ \ \ \ \ \ \ \ \ \ \ \ \ \ \leq Q-1 
\end{equation}
\end{proposition}

Here, we will fix $Q$ sufficiently close to 1 for the rest of this paper. We shall call the above coordinates $C^{k,\alpha}$-harmonic coordinates. 
Now, given $C>0, i_{0} > 0$, we can define $(k, \alpha)$-H\"older norms on the class of manifolds with $\|\nabla^{k-1} 
ric\|_{L^{\infty}} < C$ and $i(M)>i_{0}$. Start by choosing a locally finite collection of balls $B(x_{i}, \rho_{0})$ 
with $(k, \alpha)$-harmonic coordinates which cover $M$ such that the balls $B\left(x_{i}, \frac{\rho_{0}}{4}\right)$ 
are disjoint. (This is possible because the Ricci curvature is bounded.) Then define the $(k,\alpha)$-H\"older norm of 
$h \in S^{2}(M)$
\begin{eqnarray}
\|h\|_{k, \alpha}  = & \sup_{x_{i}} \left\{\sum_{1 \leq |\beta| \leq k} \rho_{0}^{|\beta|}\sup_{y \in B(x_{i}, 
\rho_{0})}|\partial^{\beta}h_{ij}(y)| \right. \\
& + \left.\sum_{|\beta | = k}
\rho_{0}^{k+|\alpha|}\sup_{y_{1},y_{2} \in B(x_{i}, \rho_{0}) }\frac{|\partial^{\beta}h_{ij}(y_{1}) - 
\partial^{\beta}h_{ij}(y_{2})|}{|y_1 -
y_2|^{\alpha}}\right\}
\end{eqnarray}
where the supremum is taken over all the balls $B(x_{i}, \rho_{0})$. We can then define the $C^{k, \alpha}$ topology 
on the space of metrics near $g_{0}$ by setting the norm of a metric $g$ near $g_{0}$ to be
\begin{equation}
\|g -g_{0}\|_{k, \alpha}.
\end{equation}

From here on out $i_{0}$, $C$, $k>2$ and $\alpha$ will all be fixed, and our reference metrics will be our approximate solutions $(M_{\sigma}, g_{\sigma})$.

Now, it is clear that on open bounded sets $\Omega \subset M$, we will have that the inclusion of $C^{k, 
\alpha}(\Omega)$ into $C^{k, \alpha^{\prime}}(\Omega)$ is compact for $\alpha^{'} < \alpha$. Furthermore, by our 
control of the metric in these coordinates we will have interior Schauder estimates for elliptic operators on bounded sets $\Omega \subset M$. (c.f. \cite{16}.)

Our analysis will take place on manifolds approximating the hyperbolic manifold whose cusps we will be closing. Thus 
we need to take into account two types of noncompact behavior: moving down the cusps makes the injectivity radius 
arbitrarily small and the diameter of the end tend toward infinity, and moving out into the expanding end makes the volume tend toward infinity. We will deal with these 
issues separately. In all of these cases, we will need infinitely many coordinate charts. Thus, to be able to get uniform control 
over the entire manifold, we shall need to have coordinate charts which are ``uniformly similar'' in some sense; i.e. 
they must be defined on balls of approximately the same size and local geometry. We will also need weight functions to deal with the infinite volume of the expanding ends, and the unbounded diameter of the filled ends.

The problem that we will encounter from the injectivity radius tending toward zero is that our coordinate charts will 
have to be made arbitrarily small as we move down the cusp. On the other hand, since the geometry is hyperbolic(or, as 
we shall see, very close to it,) we can lift to a large enough cover, and then calculate the norm on this cover. 

Thus, we take Anderson's definition from \cite{1}: 
\begin{definition}
We say that a manifold has uniformly bounded local covering geometry if, given some fixed constant $i_{0}>0$, any
ball $B(x, i_{0})$ has a finite cover $\bar{B}(\bar{x}, i_{0})$ with diameter less than 1 and $i(\bar{x}) \geq
i_{0}$.
\end{definition}

Then, define the modified $C^{k, \alpha}(M)$ norm $\tilde{C}^{k, \alpha}(M)$ to be the $C^{k, \alpha}(M)$ norm, with 
the norm being evaluated in $(k, \alpha)$-harmonic coordinates on a large enough cover if the injectivity radius is 
less than $i_{0}$. 

Now, let us define the main space which we will be working with. Let $S^{2}(M)$ be space of smooth symmetric bilinear forms on $TM$. 
\begin{definition}
Define $S^{k, \alpha}$ to be the completion of $S^{2}(M)$ with respect to the $\tilde{C}^{k, \alpha}(M)$ norm.
\end{definition}

Even though the space $S^{k, \alpha}$ is well-defined for our noncompact manifolds, it is too large for our purposes, 
since it includes many forms over whose asymptotic behavior we have very little control. Furthermore, we do not want 
to change the conformal infinity of our approximate solution when we perturb it to an exact
solution, so we want our perturbation to vanish at infinity. Both of these considerations 
lead us to the following 
definition:

\begin{definition} 
Let $\rho$ be a geodesic defining function, and let $r(x)= \log\left(\frac{2}{\rho}\right)$. For $\delta>0$, we define 
the $\delta$-weighted H\"older space $S^{k, \alpha}_{\delta}(M)$ to be 
\begin{equation}
\left\{u=e^{-\delta r}u_{0} \left| u_{0} \in S^{k, 
\alpha}(M)\right.\right\}
\end{equation}
If $u \in S^{k,\alpha}_{\delta}(M)$, define $\|u\|_{k, \alpha, \delta} = \|u_{0}\|_{k, 
\alpha}=\|e^{\delta r}u\|_{k, 
\alpha}$.
\end{definition}

Note that although this norm depends on our choice of $\rho$, the space $S^{k, \alpha}_{\delta}$ does not. $\rho$ is a 
geodesic defining function iff $|\bar{\nabla}\rho| \equiv 1$ in a neighborhood of $\partial M$. Such functions have 
the property that if $r=\log\left(\frac{2}{\rho}\right)$, then 
\begin{equation}
|\nabla r| = |\bar{\nabla}\rho| \equiv 1
\end{equation}
in some neighborhood of $\partial M$. Thus $r$ is a distance function outside some compact set. It always possible to 
construct such defining functions for AH metrics(c.f. \cite{15}.)

We then have the following analytic results

\begin{theorem} (\cite{11})
Let $\delta^{\prime} < \delta$ and $\alpha^{\prime} < \alpha$. Then the inclusion $S^{k, \alpha}_{\delta}(M) 
\rightarrow S^{k, \alpha^{\prime}}_{\delta^{\prime}}(M)$ is compact.
\end{theorem}

We will be using Bochner-technique arguments, so we shall need to use forms which are square-integrable. The following
lemma describes under which conditions this occurs.

\begin{lemma} (\cite{11}) \label{L2} Let $(M, g)$ be asymptotically hyperbolic. If $\delta 
>
\frac{n-1}{2}$ then $S^{k, \alpha}_{\delta}(M) \subset L^{2}(S^{2}(M))$.
\end{lemma}

We have yet to deal with the fact that the diameter of the filling of 
each end is tending to infinity. We will once again modify our function 
spaces, but this time the reason for the modification will not be as clear 
as the two previous ones were. The reason for this is that the change we 
are making is specifically adapted to our particular problem. Because of 
this, we will wait until we actually use it to explain the motivation. The 
basic idea is that the diameter of each of the various filled-in ends is 
determined locally on the end, and has nothing to do with those of the  
other ends. Thus, to obtain uniform bounds, we must weigh our function 
spaces with a different factor for each end.

To begin, consider any of our filled manifolds $M_{\sigma}$. For now, say we are in the filling region of the first cusp, so the $r$ coordinate lies in the interval $[r_{+}, R^{1}]$, where $R^{1}=\frac{L(\sigma^{1})}{\beta}$. Let $r_{c} \in[r_{+}, R^{1}]$. Define the function
\begin{equation}
\psi_{c} = \left\{\begin{array}{cr} 
\frac{r_{c}}{R^{1}} & \mbox{if $r_{+} \leq r \leq r_{c}$} \\
\frac{r}{R^{1}} & \mbox{if $r_{c} \leq r \leq R^{1}$} \end{array}
\right. 
\end{equation}

Assuming $N$ has $l$ cusps, we may define such a function on each filled 
cusp. Let $r_{c}=(r^{j}_{c})^{l}_{l=1}$ be a $l$-tuple such that 
$r_{c}^{j}\in[r_{+},R^{j}]$. Then define $\phi_{c}$ to be equal to the 
corresponding $\psi_{c}$ on each filling region, and equal to 1 on the rest of the manifold, up to a smoothing at each gluing torus in such a way that its $\tilde{C}^{k,\alpha}$ norm is uniformly bounded. Note that no matter what our choices are, $\phi_{c}$ is bounded above by 1, but it has no uniform positive lower bound. Note also that $\phi_{c}$ is constant near each core tori. 

\begin{definition} Let $k$, $\alpha$ and $\delta$ be as above. Consider the space $S^{k, \alpha}_{\delta}(M_{\sigma})$. Let $r_{c}=(r_{c}^{j})^{l}_{j=1}$ satisfy $r_{c}^{j}\in [r_{+}, R^{j}]$ where $R^{j}=\frac{L(\sigma^{j})}{\beta}$ is the $r$-coordinate of the $j$th gluing torus in the coordinate system of the filling manifold. We define the modified $S^{k,\alpha}_{\delta, r_{c}}$ norm of $h \in S^{k, \alpha}_{\delta}(M_{\sigma})$ to be 
\begin{equation}
\|h\|_{k,\alpha,\delta,r_{c}} = \|\phi_{c}^{-1}h\|_{k,\alpha,\delta}
\end{equation}
\end{definition}

Note that the spaces $S^{k,\alpha}_{\delta, r_{c}}$ and $S^{k, \alpha}_{\delta}$ are identical for every fixed 
$M_{\sigma}$, and that a uniform bound on the $S^{k,\alpha}_{\delta, r_{c}}$ norm of a sequence implies one on the  
$S^{k, \alpha}_{\delta}$ norm, regardless of what $r_{c}$ is. $r_{c}$ will be fixed below, in section 4.

In analogy to the $C^{k, \alpha}$ norm, we can define the $C^{k, \alpha}_{\delta, r_{c}}$ norm on metrics near our approximate solution $g_{\sigma}$ by setting 
the norm of $g$ to be
\begin{equation}
\|g-g_{\sigma}\|_{k, \alpha, \delta, r_{c}}
\end{equation}
Note that any metrics within a finite distance of each other with respect to this norm must have the same conformal infinity.

We will now discuss the operator which we shall use. To eliminate trivial kernel arising from rescalings and 
diffeomorphisms, we will use the following operator instead of the trace-free Ricci curvature $z$.

\begin{eqnarray}
\Phi_{g_{\sigma}}:S^{k,\alpha}_{\delta, r_{c}}(M_{\sigma})\longrightarrow S^{k,\alpha}_{\delta, r_{c}}(M_{\sigma}) \\
g \longmapsto ric_{g}+(n-1)g +\delta_{g}^{*}B_{g_{\sigma}}(g)
\end{eqnarray}
where $g_{\sigma}$ is the approximate solution and $B_{g_{\sigma}}$ is its associated 
Bianchi operator
\begin{eqnarray}
B_{g_{\sigma}}:  S^{2}(M) \longrightarrow \Omega^{1}(M_{\sigma}) \\
B_{g_{\sigma}}(h) = \delta_{g_{\sigma}}h + \frac{1}{2}d tr_{g_{\sigma}}h
\end{eqnarray}

It is straightforward to check that this is indeed a well-defined map between these spaces. Now the linearization of 
$\Phi$ at $g_{\sigma}$ is
\begin{equation}
D_{g_{\sigma}}\Phi(h) = \frac{1}{2}\Delta_{L}h + (n-1)h,
\end{equation}
(c.f. \cite{14},) where is the Lichnerowicz Laplacian $\Delta_{L}$ is defined as
\begin{eqnarray}
\Delta_L:  S^{2}(M) \longrightarrow S^{2}(M) \\
\Delta_{L}h = D^{*}Dh + ric \circ h + h \circ ric - 2 R(h),
\end{eqnarray}
and $R$ is the action of the curvature tensor on $S^{2}(M):$
\begin{eqnarray}
Rh(X, Y) & = & tr(\left((W,Z) \mapsto h(R(X,W)Y,Z)\right) \\
& = & \sum_{i=1}^{n}h(R(X,e_{i})Y,e_{i})
\end{eqnarray}
where $\{e_{i}\}$ is an orthonormal basis with respect to the metric at which we are linearizing. Composition of
symmetric bilinear forms is defined by using the metric to identify them with elements of Hom$(TM, TM)$.

$D_{g_{\sigma}}\Phi$ is clearly elliptic, and so $\Phi$ is elliptic near $g_{\sigma}$. The relation between $\Phi_{g_{\sigma}}$ and Einstein metrics is described by the following lemma:

\begin{lemma} (\cite{14}) \label{Biq} Let $(M, g)$ be AH, with $ric_{g}<0$. If $\Phi_{g_{\sigma}}(g)=0$ and $\lim_{r \rightarrow \infty}
\left|B_{g_{\sigma}}(g)\right| = 0$, then $ric_{g}=-(n-1)g$. In other words, $g$ is Einstein with scalar curvature 
$-n(n-1)$.
\end{lemma}

The reader may be concerned that we are only 
defining our operator $\Phi$ and our function spaces near some base metric $g_{\sigma}$. This is not an issue, since we are 
using a perturbation argument, and thus will only be working in a neighborhood of the metric we wish to perturb.

We are now ready to show why the $M_{\sigma}$ of the previous section can be called approximate solutions.
	
We can see explicitly that 
\begin{eqnarray}
\|\Phi_{g_{\sigma}}(g_{\sigma})\|_{k-2,\alpha,\delta, r_{c}} & = & \|ric _{g_{\sigma}} +(n-1)g_{\sigma}\|_{k-2,\alpha,\delta, r_{c}} 
\\
& = & C\left\|\chi^{\prime \prime} + \frac{\chi^{\prime}}{r} +  \frac{\chi}{r^{2}}\right\|_{k-2,\alpha,\delta} \\
&\leq & C_{1}\frac{1}{|\sigma|^{n-1}}
\end{eqnarray}

Note that one must be careful in the last line, since we are bounding H\"older norms, in which derivatives are 
calculated with respect to harmonic coordinates, and not with respect to the coordinate $r$. But $r$ is related to the 
geodesic coordinate $s$ by 
\begin{equation}
r = O(e^{s})
\end{equation}
for large $s$, so we can establish the bound with respect to $s$, and then translate back into terms of $R$.

	Thus, we have the following proposition:
\begin{proposition}
Let $(N,g)$ be a geometrically finite hyperbolic manifold, whose $l$ cusps all have toric cross-sections. Then for any 
$l$-tuple $\sigma$ of geodesics, with the $i$th geodesic drawn from the $i$th cusp cross-section, it is 
possible to 
construct a manifold $(M_\sigma, 
g_{\sigma})$ such that 
\begin{equation}
\|\Phi_{g_{\sigma}}(g_{\sigma})\|_{k-2,\alpha,\delta, r_{c}} = O(|\sigma|^{1-n})
\end{equation} 
These $(M_\sigma, g_{\sigma})$ are AH and have the same conformal infinity as $(N,g)$.
\end{proposition}

We will also need elliptic estimates for these operators. In analogy to the weighted H\"older norms, define the 
$(\delta, r_{c})$-weighted $L^{\infty}$ norm to be $\|h\|_{L^{\infty}_{\delta,r_{c}}} = \|e^{\delta r}\phi_{c}^{-1}h\|_{L^{\infty}}$.
\begin{proposition} 
Let $L_{g} = D_{g}\Phi_{g_{\sigma}}$ be the linearization of $\Phi_{g_{\sigma}}$ at $g$, where $\|g-g_{\sigma}\|_{k, \alpha, 
\delta, r_{c}}<\epsilon_{0}$, and $\epsilon_{0}>0$ is chosen such that $L_{g}$ is elliptic. Then there is some constant 
$\Lambda_{0}$, depending on $k, \alpha, \epsilon_{0}$ and $\delta$ such that we have the following estimate
\begin{equation}
\|h\|_{k, \alpha, \delta, r_{c}} \leq \Lambda_{0}(\|Lh\|_{k-2, \alpha, \delta,r_{c}} + \|h\|_{L^{\infty}_{\delta, r_{c}}})
\end{equation}
\end{proposition}

The proof of this is essentially that of the analogous estimate without the $r_{c}$'s in \cite{11}. The critical property of the weight functions $e^{-\delta r}$ and $\phi^{-1}_{c}$ that is required is that one can establish an estimate of the form
\begin{equation}
\|f\|_{C^{2}} \leq C\|f\|_{C^{0}}
\end{equation}
with $C$ independent of the manifold and the choice of the $r_{c}$. Then, we can apply the standard Schauder estimate to the weighted functions, and the above estimate allows us to interchange the weight functions and the operator $L$. The point to note here is that $\Lambda_{0}$ does not depend in any way on $\sigma$ or $r_{c}$.

We end this section by quoting a key result of Graham and Lee's(\cite{15}) on the behavior of the operator 
$L=D_{g_{\sigma}}\Phi_{g_{\sigma}}$. We 
know that $L$ is elliptic, but since we are working on noncompact manifolds $M$, it does not follow that $L$ is 
Fredholm. Thus if we do not choose our function spaces carefully, $\ker(L)$ or coker$(L)$ could well be 
infinite-dimensional. Biquard's results give appropriate conditions to guarantee that $L$ is Fredholm. (Also see 
\cite{14}).

\begin{theorem} (\cite{15}) \label{GL}
Let $(M^n, g_{0})$ be an asymptotically hyperbolic manifold. If $\delta \in (0,n-1)$, then
\begin{equation} 
L: S^{k, \alpha}_{\delta}(M^{n}) \longrightarrow S^{k-2, \alpha}_{\delta}(M^{n})
\end{equation}
is Fredholm. Furthermore, $L$ is an isomorphism iff $\ker_{L^2}(L)=0$.
\end{theorem}

We will need to use both Lemma \ref{L2} and Theorem \ref{GL} below, so we will restrict 
$\delta$ to the interval $\left(\frac{n-1}{2},n-1\right)$.

\section{Control of Inverse}

In this section, we will perform the analysis which will allow us to invert our operators $\Phi_{g_{\sigma}}$. To do this, we will need to invert the linear operators $L_{\sigma} = 2D_{g_{\sigma}}\Phi_{g_{\sigma}}$ and get some kind of uniform control on the behavior of the inverses. We cannot get an absolute uniform bound on their operator
norms, but we can make sure that they do not grow too fast with respect to the length of the Dehn surgery $\sigma$. We prove via a contradiction argument that we may choose $r_{c}$ such there exists a $\Lambda$ independent of $\sigma$ such that
\begin{equation} \|h\|_{k, \alpha,\delta, r_{c}} \leq \Lambda \|L_{\sigma}(h)\|_{k-2,\alpha,\delta, r_{c}} \label{est}
\end{equation}
for all $h \in S^{k, \alpha}_{\delta, r_{c}}$, provided $k\geq 2$ and that $\sigma$ is large enough

The above estimate shows that $\ker(L_{\sigma})=0$, so by theorem \ref{GL} we get that
\begin{equation}
L_{\sigma}^{-1}:S^{k-2,\alpha}_{\delta, r_{c}}\longrightarrow S^{k,\alpha}_{\delta, r_{c}}
\end{equation}
 is well-defined and
\begin{equation}
\|L^{-1}_{\sigma}f\|_{k, \alpha, \delta, r_{c}} \leq \Lambda  \|f\|_{k-2, \alpha, \delta, r_{c}} 
\end{equation}

The proof of estimate (\ref{est}) is essentially identical to Anderson's proof in the finite-volume case(\cite{1}). The only additional
difficulty that we face here is that $S^{k,\alpha}_{\delta, r_{c}}$ functions are not 
necessarily in $L^{2}$, but as mentioned above, by choosing $\delta >
\frac{n-1}{2}$ we can avoid this problem. In this section, any unlabeled norms will be assumed to be $L^{2}$ norms and 
$\delta \in \left(\frac{n-1}{2},n-1\right)$ will be assumed to be fixed.

\begin{proposition} Let
$(M_{\sigma}, g_{\sigma})$ be a sequence of approximate solutions. Then there exists some $r_{c}$ and a constant $\Lambda$ independent of
$\sigma$ such that

\begin{equation}
\|h\|_{k, \alpha, \delta, r_{c}} \leq \Lambda \|L_{\sigma}(h)\|_{k-2, \alpha, \delta, r_{c}}
\end{equation}

\end{proposition}

\textbf{Proof:} We work by contradiction, so we will have to
take limits. This leads to some difficulty, since there is no uniform bound on the diameter of the $M_{\sigma}$'s, and so the limits will not be uniquely defined. On the other hand, all of the limits are Einstein, which gives them extra structure which we will exploit.

Let us set up the contradiction.  
If there is no $\Lambda$ such that 
\begin{equation}
\|h\|_{k, \alpha,\delta, r_{c}} \leq \Lambda \|L_{\sigma}(h)\|_{k-2, \alpha,\delta, r_{c}},
\end{equation}
for all $\sigma$, then necessarily there is a sequence $h_{i} \in S_{\delta, r_{c}}^{k, \alpha}(M_{\sigma_{i}}, g_{\sigma_{i}})$ such that 
\begin{equation} 
\|h_{i}\|_{k, \alpha,\delta, r_{c}} = 1
\end{equation}                                                                  but
\begin{equation}
\|L_{i}(h_{i})\|_{k-2, \alpha,\delta, r_{c}} \longrightarrow 0
\end{equation}
where we have replaced the subscript $\sigma_{i}$ by an $i$.

By the Schauder estimates, we know that there is a $\Lambda_{0}$ independent of $\sigma$ and $r_{c}$ such that
\begin{equation} 
\|h\|_{k, \alpha, \delta, r_{c}} \leq \Lambda_{0} (\|L_{\sigma}(h)\|_{k-2, \alpha, \delta, r_{c}} +
\|h\|_{L^{\infty}_{\delta, r_{c}}})                                   
\end{equation}      
for all $h \in S^{k, \alpha}_{\delta, r_{c}}(M_{\sigma})$. Thus, under our assumption,
\begin{equation}
\|h_{i}\|_{L^{\infty}_{\delta, r_{c}}}\geq \Lambda^{-1}_{0} > 0
\end{equation}
Therefore, showing that $\|h_{i}\|_{L^{\infty}_{\delta, r_{c}}} \rightarrow 0$ will give us a contradiction.

We will suppress the $r_{c}$ for the next while, since it will not play a role until later.

The most natural way of looking at $M_{\sigma}$ is to see it as being made up of two distinct pieces: the original hyperbolic manifold $N$ and the collection of black hole metrics with which we are filling in the cusps of $N$. Our strategy is to take the limit of the $h_{i}$'s, which will lead to infinitesimal Einstein deformations of each piece. We will get our contradiction by showing that there can be no nontrivial deformations. 

We will spend most of our time working on the filled regions of the cusps, i.e. the ends which are close to the black hole metrics. We will use the variables $r$ to refer to the $r$-variable in our parametrization of the black hole metrics, and $R$ will refer to the gluing region, as seen from the black hole metrics. Note that, by construction, $R=\frac{L(\sigma)}{\beta_{1}}$. Our analysis will take place on each black hole region separately, so there will be no risk of confusion. 

For the region associated to the $k$th cusp, we have the relation:
\begin{equation}
R^{k}_{i} = \beta_{1}L(\sigma^{k}_{i}) \geq \beta_{1} |\sigma_{i}|
\end{equation}

To begin, note that we have the following Weitzenb\"ock formula(c.f. \cite{4}):

\begin{equation}
\delta d h + d \delta h = D^* D h - R h + h \circ ric
\end{equation}
where $d$ is the exterior derivative on vector-valued one-forms induced by the connection, and $\delta$ is its formal adjoint.

From this, we get that 
\begin{eqnarray}
Lh &=  & D^{*}Dh - 2Rh +h \circ ric  +  ric \circ h + 2(n-1)h \\
& = & \delta d h + d \delta h - R h + ric \circ h + 2(n-1)h
\end{eqnarray}

We will work primarily with this form of $L$, since we can prove stronger positivity properties with it.

Now, by construction, we have that on $M_{i}$, 

\begin{equation}
ric + (n-1)g = \tau(r)
\end{equation}
where $\tau$ is supported on the region of the black hole metric $(R_{i} - 2, R_{i}-1)$ and satisfies
\begin{equation}
|\tau| \leq C R_{i}^{1-n}
\end{equation}

Thus, 
\begin{equation}
Lh = \delta d h + d \delta h - R h + (n-1)h + \tau \circ h
\end{equation}
 
Now, consider
\begin{equation}
\langle Lh, h \rangle = \langle \delta dh, h \rangle + \langle d\delta h, h \rangle - \langle Rh, h \rangle + (n-1)\|h\|^{2} + \langle \tau \circ h,h\rangle
\end{equation}

Recall that we have $Lh \rightarrow 0$ and $h$ bounded, and we want to show that $h \rightarrow 0$. This should be possible as long as all the terms on the right hand side 
are positive or tend toward $0$. Integration by parts will work on the first two terms, so we need to get a handle on the term $\langle Rh, h \rangle$. We will do this by controlling the pointwise norm $(Rh, h)$. First, we will break it up 
into three pieces. Let $h = h_{0} + \frac{(tr h)g}{n}$ (so $h_{0}$ is the trace-free part of $h$.)

\begin{lemma} \label{decomp}
\begin{equation}
(Rh, h) = (Rh_{0}, h_{0}) + \mu_{i}(r) + O\left((tr h)^{2}\right) 
\end{equation}
Where $\mu_{i}(r)$ is supported in the black hole region and is $O\left(R_{i}^{-(n-1)}\right)$
\end{lemma}
\textbf{Proof:} To begin, note that 
\begin{eqnarray}
(Rh, h) & = & \left(R\left(h_{0} + \frac{(tr h)}{n}g\right),  h_{0} + \frac{(tr h)g}{n}\right) \\
	& = & (Rh_{0}, h_{0}) + \frac{tr h}{n}(Rg, h_{0}) + \frac{tr 
h}{n}(Rh_{0}, g)\\
& &  + \frac{(tr h)^{2}}{n^{2}}(ric, g)
\end{eqnarray}

Now, we have that $(ric,g)=s$, which is uniformly bounded. Furthermore, if we take an orthonormal frame $(e_{j})$ which diagonalizes the curvature tensor, 
\begin{equation}
Rh(e_{j}, e_{k}) = \sum_{l}h(R(e_{j}, e_{l})e_{k}, e_{l})
\end{equation}
Thus,
\begin{eqnarray}
(Rg, h_{0}) & = & (g, Rh_{0}) \\
  & = & tr(Rh_{0}) \\
  & = & \sum_{j,k}h_{0}(R(e_{j}, e_{k})e_{j}, e_{k}) \\
  & = & \sum_{j,k}K_{jk}h_{0}(e_{k},e_{k}) \\
  & = & -(n-1) tr(h_{0}) +  \mu_{i} \\
  & = & \mu_{i}
\end{eqnarray}
Where $\mu_{i}(r) = \sum_{j,k}(\delta_{j}^{k}+K_{jk})h_{0}(e_{k},e_{k})= O\left(R_{i}^{-(n-1)}\right)$ since $h_{0}$ is 
uniformly bounded.
\qed

The following lemma allows us to control $(Rh_{0},h_{0})$.

\begin{lemma} \cite{4} \label{eigen} Let 
\begin{equation}
a = \sup_{\{h_{0} | tr h_{0} = 0\}}  
\frac{(Rh_{0},h_{0})}{|h_{0}|^{2}}
\end{equation} 
be the largest eigenvalue of $R$ acting on $S^{2}_{0}$. Then 
\begin{equation}
a < (n-2)K_{max} - ric_{min}
\end{equation}
\end{lemma}

For the next term, we can use the following result from \cite{1}.

\begin{lemma}
As $i \rightarrow \infty$, we have $\|tr h_{i}\|_{L^{2}} \longrightarrow 0$
\end{lemma}

Now, let $\mathcal{U}_{\rho} = \{x|r(x) < \rho\}$ be a tubular neighborhood of 
the totally geodesic core $T^{n-1}$'s, and let $M^{\rho}_{i} = M_{i} - \mathcal{U}_{\rho}$. We will fix $\rho\leq 
R_{i}=\inf_{k} R_{i}^{k}$ 
below. Note that $\mathcal{U}_{\rho}$ has $q$ connected components, where 
$q$ is the number of cusps of $N$.

By Lemmas \ref{decomp} and \ref{eigen},  on $M_{i}^{\rho}$.
\begin{eqnarray}
(R(h), h) & = & (R(h_{0}), h_{0}) + \mu_{i} + O\left((tr h)^{2}\right) \\
	& \leq & \left((n-2)K_{max} - ric_{min}\right)|h|^{2} + \mu_{i} + C_{1}|tr h|^{2} \\
	& = & \left(-(n-2)+ (n-1)+ O(R_{i}^{1-n})\right)|h|^{2} \\
& & \ \ + \mu_{i} + C_{1}|tr h|^2\\
&\leq& \left(1 +  C_{2}R_{i}^{1-n}\right)|h|^{2}  +  \mu_{i} + C_{1}|tr 
h|^2 \label{bound}
\end{eqnarray}

Now consider
\begin{eqnarray}
\int_{M^{\rho}} (Lh, h) dV & = & \int_{M^{\rho}} (\delta d h, h) + (d \delta h, h)  - (R h, h)  \\ 
& & + (n-1)|h|^{2}+(\tau \circ h, h) \  dV \\
	& = & \int_{M^{\rho}} |\delta h|^{2} + |dh|^{2} - (R h, h) + 
(n-1)|h|^{2}  \\ 
 & & + (\tau\circ h,h) \  dV 
+ \int_{\partial \mathcal{U}_{\rho}} Q(h, \partial h) dA 
\end{eqnarray}

	Here, $Q(h, \partial h)$ is the boundary term from the integration by parts. It is a fixed quadratic polynomial 
in $h$ and its derivatives. Choose some $\epsilon>0$. By our estimates on $(Rh, h)$, and assuming that $i$ is large enough that $|\tau_{i}|<\epsilon$ and
\begin{equation}
1 + C_{2}R_{i}^{1-n} \leq \frac{5n}{12}
\end{equation}
we get that this last quantity is 
\begin{eqnarray}
& \geq &  \int_{M^{\rho}} - \left(\frac{5n}{12} + \epsilon\right)|h|^{2} + 
(n-1)|h|^{2}\  dV -
C_{1}\int_{M^{\rho}}(trh)^{2}\ dV \\ 
& & \ \ - \int_{M^{\rho}} \mu_{i} \ dV + 
\int_{\partial \mathcal{U}_{\rho}} Q(h, \partial h)\ dA
\end{eqnarray}

Now, $\mu_{i}$ is $O(R_{i}^{1-n})$ and supported on a region of bounded 
volume(the black hole region) and by the 
previous lemma $\|trh_{i}\|_{L^{2}}\longrightarrow 0$, so we may choose 
$i$ large enough that the previous quantity is
\begin{eqnarray}
& \geq &  \int_{M^{\rho}} - \left(\frac{5n}{12} + \epsilon\right)|h|^{2} + 
(n-1)|h|^{2}\  dV  - 2\epsilon \\ 
& & \ \ \ \ + 
\int_{\partial \mathcal{U}_{\rho}} Q(h, \partial h) dA \\ 
& = & \left(\frac{7n}{12} -1 -\epsilon\right)\int_{M^{\rho}} |h|^{2} \ dV  + 
\int_{\partial \mathcal{U}_{\rho}} Q(h, \partial h) dA - 2\epsilon
\end{eqnarray}
We also have that 
\begin{equation}
\int_{M^{\rho}} (Lh, h) dV \leq \frac{1}{2}\left(\int_{M^{\rho}} |Lh|^{2} dV + \int_{M^{\rho}} |h|^{2} dV\right)
\end{equation}
and 
\begin{equation}
\left|\int_{\partial \mathcal{U}_{\rho}} Q(h, \partial h) dA\right| \leq C Vol(\partial \mathcal{U}_{\rho})
\end{equation}
 since we have $C^{k, \alpha}$ control over the $h_{i}$'s.

Combining all this information gives
\begin{equation}
\frac{1}{2}\int_{M^{\rho}} |Lh|^{2}\  dV + C Vol(\partial \mathcal{U}_{\rho}) \geq 
\left(\frac{7n}{12}-\epsilon-\frac{3}{2}\right) \int_{M^{\rho}} |h|^{2}\  
dV -2\epsilon
\end{equation}
where we can make $\epsilon>0$ arbitrarily small by making $i$ large.

By assumption, $\|Lh\|_{L^{\infty}_{\delta}} \longrightarrow 0$, so we 
have that $\int_{M^{\rho}} |Lh|^{2} dV \longrightarrow 0$. Remark that $inj(\partial \mathcal{U}_{\rho})= 
O\left(\frac{\rho}{R_{i}}\right)$. Thus, if we choose a sequence $\rho_{i}$ such that 
$\frac{\rho_{i}}{R_{i}}\longrightarrow 0$, we'll obtain that $\int_{M^{\rho_{i}}} |h|^{2} dV$ tends to $0$. Since 
$S^{k, \alpha}_{\delta} \subset L^{2}$, this gives us uniform convergence of the $h_{i}$'s on any set whose 
injectivity radius remains bounded below. Thus, by a diagonal argument, we can find some sequence $r_{i}$ such that 
$h_{i}$ converges uniformly to 0 on $M^{r_{i}}$ and $\frac{r_{i}}{R_{i}} \longrightarrow 0$. Note that we are finding a different $r_{i}$ for each cusp.

We will use these $r_{i}$ as the parameter for the cusp we are on in the weight function $\phi_{c}$. In other words, the $r_{i}$ obtained for the $j$th cusp of the manifold $M_{i}$ will be the $r^{j}_{c}$ which determine the weight function on that cusp.

	Now, let us examine what is happening on each component of $\mathcal{U}_{r_{i}}$, the complement of the set 
$M^{r_{i}}$. By construction, the core tori are collapsing to points, so any neighborhood of these tori is degenerating 
to a line segment. Since we want to have a nice limit, we lift everything to finite covers in order to ``unwrap the 
collapse.''(c.f. \cite{2}.) Choose a sequence of points $p_{i}$ in a core torus. Since $inj(T(r))= 
O\left(\frac{r}{R_{i}}\right)$, by lifting $(\mathcal{U}_{r_{i}}, g_{i}, p_{i})$ to an 
$\left[\frac{R_{i}}{r_{i}}\right]$-fold cover, where $[\ \ ]$ is the greatest integer part function, we will get a 
sequence of manifolds whose injectivity radius is bounded away from $0$. 

	By definition of the $S^{k, \alpha}_{\delta, r_{c}}$ norm, the $C^{k, \alpha}$-norm of these manifolds is 
bounded, so we 
get a limit manifold $(M_{BH}, g_{BH})$. Clearly the $h_{i}$'s lift to the finite 
covers too, so we get lifted forms $\tilde{h_{i}}$. These $\tilde{h_{i}}$ satisfy \begin{equation}
 \|\tilde{h_{i}}\|_{k, \alpha, \delta} \leq C
\end{equation}
on the lifted manifolds $(\tilde{\mathcal{U}}_{r_{i}},\tilde{g_{i}})$.   
Thus, given $\alpha^{\prime}<\alpha$, and 
$\delta^{\prime} < 
\delta$, we can extract a subsequence to get a limit $\tilde{h} \in S^{k,\alpha^{\prime}}_{\delta^{\prime}}$. Note that $\tilde{h}$ must be $T^{n-1}$ invariant and satisfy $L\tilde{h}=0$. 

If $|r_{i}-r_{+}|$ is uniformly bounded, the torus $T(r_{i})$ stays within a fixed distance of the core torus for all 
$i$. Thus, the core torus can always see the region on which $h_{i}$ is tending uniformly to 0. We can therefore take a 
 pointed limit based at $p_{i} \in T(r_{i})$, and conclude that $h=0$ if $r>r_{i}$. By analyticity of infinitesimal 
Einstein deformations, this leads us to conclude that $h$ is identically $0$, giving us our contradiction. Thus, we 
shall assume that $r_{i} \longrightarrow \infty$. This gives us that the limiting manifolds $(M_{BH}, g_{BH})$ are complete. We are going to work at an infinite distance from the conformal infinity of $M_{\sigma}$, so we can drop the weight factor $\delta$. It will also be understood that we are working on the lifted manifolds, so we will suppress the tildes.

At this point, we would like to say that since $h_{i}(r_{i}) \longrightarrow 0$ and $r_{i} \longrightarrow \infty$, this forces our limiting sequence to have 
\begin{equation}
\lim_{r \rightarrow \infty} \|h\| = 0
\end{equation}
We could then apply the following results to get our contradiction:
\begin{proposition}
$h$ is tangent to the space of $T^{n-1}$-invariant AHE metrics on $M$
\end{proposition}
\textbf{Proof:} This is nontrivial, since spaces of AHE metrics are 
infinite-dimensional, and so vector fields do not necessarily integrate. 
By \cite{3}, however, we know that if the AHE manifold $(M,g)$ has a $C^{2}$ conformal compactification, then infinitesimal deformations do indeed integrate. The function $\rho=r^{-1}$ is certainly $0$ exactly at the conformal infinity of $(M_{BH},g_{BH})$, and with respect to the compactified metric,
\begin{equation}
|d\rho|_{\bar{g}}^{2} = \left|\frac{dr}{r^{2}}\right|_{\bar{g}}^{2} = \frac{V(r)}{r^{2}} 
\end{equation}
is nonzero on the boundary, so $\rho$ is a defining function.  Near the boundary, the compactified metric is
\begin{equation}
\bar{g} = \frac{dr^{2}}{r^{2}V(r)} + \frac{V(r)}{r^{2}}d\theta^{2} + g_{T^{n-2}}
\end{equation}
One may replace $r$ by the coordinate $s$, where 
\begin{equation}
\frac{ds}{dr} = \frac{1}{r\sqrt{V(r)}}
\end{equation}
To get
\begin{equation}
\bar{g} = ds^{2} + F(s)d\theta^{2} + g_{T^{n-2}}
\end{equation}
and a short calculation shows that $F$ is $C^{2}$ up to the boundary.
\qed

\begin{proposition}(\cite{17} c.f. also \cite{40}) \label{rigidity}
Let $g$ be a complete $T^{n-1}$-invariant AHE metric on the solid torus $D^{2}\times T^{n-2}$. Then $g$ is a black 
hole metric.
\end{proposition}

\begin{proposition} \label{unique} Let $g_{t}$ be a curve of complete $T^{n-1}$-invariant 
AHE metric on 
the solid torus $D^{2}\times T^{n-2}$. Then $g_{t}$ is completely determined by $g_{0}$ and the curve $\gamma_{t}$ consisting of the conformal infinities of the $g_{t}$'s.
\end{proposition}

\textbf{Proof:} By the previous proposition, all the $g_{t}$'s are covered by $(D^{2}\times\R^{n-2}, \tilde{g}_{BH})$. Thus,
\begin{equation}
g_{t} = \frac{\tilde{g}_{BH}}{\Gamma_{t}}
\end{equation}
where $\Gamma_{t}\subset Isom(g_{BH})$ is isomorphic to $\Z^{n-2}$. 
Consider the ``defining function'' $\rho = r^{-1}$ for 
$\tilde{g}_{BH}$(note that $(D^{2}\times\R^{n-2}, 
\rho^{2}\tilde{g}_{BH})$ is not a compact manifold with boundary.)

 Then the conformal infinity of  $(D^{2}\times\R^{n-2}, 
\rho^{2}\tilde{g}_{BH})$ is a flat $S^{1}\times\R^{n-2}$. The 
action of $\Gamma_{t}$ commutes with multiplication by $\rho$, so the conformal infinity of $g_{t}$ is the quotient of 
the flat $S^{1}\times\R^{n-2}$ by $\Gamma_{t}$. Conversely, the conformal infinity $\gamma_{t}$ of $g_{t}$ also 
determines the group $G_{t} = \Gamma_{t}+ v_{t}$ up to conjugacy, where 
\begin{equation}
\gamma_{t} = \frac{\R^{n-1}}{G_{t}}
\end{equation}
Given an initial $g_{0}$, we can identify $v_{0}$ with the $S^{1}$ in the universal cover. This determines $v_{t}$ for 
$t>0$, and thus $\Gamma_{t}$ determines $g_{t}.$
\qed

Combining these propositions gives us the following corollary:
\begin{corollary} \label{defs}
A $T^{n-1}$-invariant AHE metric on a solid torus has no nontrivial infinitesimal deformations for which
\begin{equation}
\lim_{r \rightarrow \infty} |h(r)| = 0 
\end{equation}
\end{corollary}

Unfortunately(for us,) we cannot apply these propositions directly; the problem comes from the fact that the $h_{i}$'s 
tend to $h$ uniformly on compact subsets, and we cannot relate the rate of convergence to the size of the compact set. 
Thus, we cannot know that $r_{i}$ is included in each set. Consider the following example: Let $I_{k}=[-k, k]$, and 
\begin{eqnarray*}
f_{k}: I_{k} \longrightarrow [-1, 1] \\
x \mapsto \frac{x}{k}
\end{eqnarray*}

	Clearly, $\|f_{k}^{\prime}(x)\|_{L^{\infty}} \longrightarrow 0$ as $k \longrightarrow \infty$, so one would expect that the $f_{k}$'s are tending toward a constant function. This constant may, however, depend on the basepoint $x_{k}$; let $x_{k}=\alpha k$, where $-1 < \alpha <1$. Then the (pointed) Gromov-Hausdorff limit of the triple $(I_{k}, x_{k}, f_{k}(x))$ will be the triple $(\R, 0, \alpha)$. Thus the limit of the $f_{k}$ is indeed a constant function, but the constant depends on the choice of the basepoints $x_{k}$. 

The issue here is that the manifolds $I_{k}$ are converging to their limit faster than the $f_{k}$ are converging to theirs, so the convergence cannot be made uniform on the whole set. 

More precisely, we have that
\begin{eqnarray}
|f_{k}(x) - f_{k}(y)| & \leq & \left|\int_{x}^{y}f^{\prime}_{k}(s) \ ds \right| \\ & \leq & \|f_{k}^{\prime}\|_{L^{\infty}}\left|\int_{x}^{y}\ ds \right| \\ 
& \leq & \|f_{k}^{\prime}\|_{L^{\infty}}|x-y|
\end{eqnarray}

Therefore, if we choose points $x_{k}, y_{k}$ whose distance is increasing fast enough, we cannot conclude that 
$f_{k}(x_{k})$ and $f_{k}(y_{k})$ have the same limit

On the other hand, if we require that $\|f_{k}^{\prime}\|_{L^{\infty}} \longrightarrow 0$ more rapidly than any two 
points can separate, say by demanding that $diam(I_{k})\|f'_{k}\|_{L^{\infty}} \longrightarrow 0$, then we can get that 
\begin{eqnarray}
|f_{k}(x) - f_{k}(y)| & \leq & \|f_{k}^{\prime}\|_{L^{\infty}}|x-y| \\
& \leq & diam(I_{k})\|f'_{k}\|_{L^{\infty}} \longrightarrow 0
\end{eqnarray}
no matter which basepoints we take. 

Let us now attempt to adapt this argument to the operator $L$. We will start by analyzing 
the limiting case to motivate our choice of the weight $\phi_{c}$, and then we will get our 
hands dirty with the actual estimates that we need to finish our proof. 
As we will see below, on 
$T^{n-1}$-invariant deformations $h$, the components of $Lh=0$ are asymptotic to Euler equations. For the components $(Lh)_{1j}$, this equation has a nontrivial 0th order term, so all of its solutions either blow up or decay to $0$. Its other components are asymptotic to equations of the form 
\begin{equation}
Lf = r^{2}f^{\prime\prime} + nrf^{\prime} = 0
\end{equation}
which have constant solutions. This leads to a problem for us, since we could have the 
same situation as above; even though $f(r_{i})\rightarrow 0$ and $f$ tends to a constant, 
we cannot conclude that $f\longrightarrow 0$ everywhere. This is where the weight $\phi_{c}$ in our norms comes in.

By use of an integrating factor(\cite{26}), we may rewrite this as 
\begin{equation}
f(r)= \int \frac{1}{r^{n}}\int s^{n-2}Lf(s) \ ds \ dr
\end{equation}

Thus,
\begin{eqnarray}
|f(r_{1})-f(r_{2})| & = & \left| \int_{r_{1}}^{r_{2}} \frac{1}{r^{n}}\int_{r_{+}}^{r} s^{n-2}Lf(s) \ ds \ dr \right| \\
& \leq & \frac{\|Lf\|_{L^{\infty}}}{n-1} \left|\int_{r_{1}}^{r_{2}}r^{-1} + O\left(r^{1-n}\right)\ dr \right| \\
& \leq & \frac{\|Lf\|_{L^{\infty}}}{n-1} \left| \log\left(\frac{r_{2}}{r_{1}}\right)+ C_{0} \right| \\
& \leq & C_{1}\frac{R_{i}}{r_{i}}\|Lf\|_{L^{\infty}} \\
& \leq & C_{1}\|\phi^{-1}_{c}Lf\|_{L^{\infty}} \longrightarrow 0
\end{eqnarray}

 As we remarked above, our situation is a bit more delicate, since we are not working 
exactly with this operator, 
but with perturbations of it. Furthermore, since we need to use the rate at which the $h_{i}$'s converge, we cannot just work with limits, but must rather get precise bounds on how things behave asymptotically.

We now want to analyze the system of ODE's that these deformations must satisfy. The proof of the following 
proposition consists of a long calculation which is quite complicated, but fairly straightforward. The interested 
reader may consult \cite{40} to see the gory details, and the masochistic reader may attempt it for her- or himself.
	
\begin{proposition} \label{deform}
Say we have a black hole metric $g$.
Let $e_{1}=\sqrt{V}\partial_r$, $e_{2}=\frac{1}{\sqrt{V}}\partial_{\phi_{j}}$, and 
$e_{j}=\frac{1}{r}\partial_{\phi_{j}}$, where the $\partial_{\phi_{j}}$'s, $3\leq j \leq 
n$ form an orthonormal basis 
for the core torus. Then if $h$ is $S^{1}\times T^{n-2}$-invariant, $Lh=2D_{g}\Phi_{g}h$ is given by 
\begin{eqnarray}
(Lh)_{11}  & = &  Ah_{11} + h_{11}\left(\frac{(V^{\prime})^{2}}{2V} + 
\frac{2(n-2)V}{r^{2}}\right)\\
& &  -h_{22}\left(\frac{(V^{\prime})^{2}}{2V} + 
2K_{12}\right)\\& &  - 2\sum_{k>2}\left(\left(\frac{V}{r^{2}}\right) + K_{1k}\right)h_{kk} \\
(Lh)_{22} & =  &  Ah_{22} +h_{22}\frac{(V^{\prime})^{2}}{2V} -h_{11}\left(\frac{(V^{\prime})^{2}}{2V}+ 2K_{12}\right) \\
 & & \ \ \ \ \ \ \ \ \ \ \ \  -2\sum_{k>2}K_{2k}h_{kk} \\
(Lh)_{12} & = & Ah_{12} + h_{12}\left(\frac{(V^{\prime})^{2}}{V} + \frac{2(n-2)V}{r^{2}}+2K_{12}\right)
\end{eqnarray}
where 
\begin{equation}
Ah_{ij} =  \left(-Vh^{\prime \prime}_{ij} - \left(V^{\prime}+  \frac{n-2}{r}V\right)h_{ij}^{\prime}\right)
\end{equation}

If $j>2$, we have 
\begin{eqnarray}
(Lh)_{jj} & = & Ah_{jj} + \frac{2V}{r^{2}}(h_{jj} - h_{11}) - 2\sum_{k\neq j}K_{kj}h_{kk} \\
(Lh)_{1j} & = & Ah_{1j} + h_{1j}\left(\frac{(V^{\prime})^{2}}{4V} + \frac{(n+1)V}{r^{2}}+2K_{1j}\right) \\
(Lh)_{2j} & = & Ah_{2j} + h_{2j}\left(\frac{(V^{\prime})^{2}}{4V} + \frac{V}{r^{2}}+2K_{2j}\right)
\end{eqnarray}
and finally, if $i,j>2$, we get
\begin{equation}
(Lh)_{ij} = Ah_{ij} + h_{ij}\left(\frac{(V^{\prime})^{2}}{2V}+2K_{ij}\right)
\end{equation}

\end{proposition}

This system seems unmanageable for the black hole metrics, but we can get around this by noting the following two facts:
\begin{proposition}
Let $g_{C}$ be a complete hyperbolic cusp metric
\begin{equation}
g_{C} = r^{-2}dr^{2} + g_{T^{n-1}}
\end{equation}
where $g_{T^{n-1}}$ is an arbitrary flat metric on the torus with orthonormal basis $\partial_{\phi_{j}}$, 
$2\leq j\leq n$. Let 
$\left(e_{1}=r\partial_{r}, e_{j}=\frac{\partial_{\phi_{j}}}{r}; 
j\geq 2\right)$ form an orthonormal frame for $g_{C}$.  
Then if 
$h$ is $T^{n-1}$-invariant, $L_{C}h=L_{g_{C}}h$ is given by the following formulae. \begin{equation}
(Lh)_{11} = Ah_{11} + 2(n-1)h_{11}
\end{equation}
and if $j,k>2$,
\begin{eqnarray}
(Lh)_{jj} & = & Ah_{jj}  + 2 tr h -2h_{11} \\ 
(Lh)_{1j} & = & Ah_{1j} + nh_{1j} \\
(Lh)_{jk} & = & Ah_{jk}
\end{eqnarray}
where $A = -r^{2}\partial_{r}^{2} - nr\partial_{r}$.
\end{proposition}
\textbf{Proof:} Set $V(r) =r^{2}$ and $K_{jk}=-1 + \delta_{j}^{k}$ above. \qed
	
Now that we have the much simpler form of $L$ for $g_{C}$, we must relate it to the corresponding operator on the black hole metrics. We will use $C^{k}$ estimates instead of H\"older ones, since they are easy to establish, and we do not need the stronger norms; we already have the existence of the limit form $h$. The $C^{k}$ norms will be calculated in the same harmonic coordinates as the $S^{k, \alpha}_{\delta}$ ones. An easy calculation gives:

\begin{proposition}
If $\|h\|_{C^{2}}$ is bounded and $T^{n-1}$-invariant, then 
\begin{equation}
\|L_{C}h - L_{BH}h\|_{L^{\infty}} = O(r^{-(n-1)})
\end{equation}
\end{proposition}

Thus, if $h$ is $T^{n-1}$-invariant and bounded in $C^{2}$, we get the following systems of equations for $Lh = 0$:
\begin{eqnarray}
Ah_{11} + 2(n-1)h_{11} =  u_{11} \\
Ah_{jj}  + 2 tr h -2h_{11} = u_{jj}\\ 
Ah_{1j} + nh_{1j} = u_{1j} \\
Ah_{ij} = u_{ij}
\end{eqnarray}
where $i,j>1$ all the components of $h$ are bounded, and $|u| = O(r^{-(n-1)})$

Recall that we want to prove that  
\begin{equation}
\lim_{r\rightarrow\infty} |h(r)| = 0
\end{equation}

This is straightforward in the case of the components $h_{11}$, $h_{1j}$ and $h_{jj}$. All of these satisfy nonhomogeneous Euler equations with nonzero indicial roots. One can write out the solutions to the nonhomogeneous equations in terms of the fundamental solutions using variations of parameters(c.f. \cite{26},) and see that bounded solutions must tend to 0 as $r\rightarrow \infty$.

Thus, all that is left is to examine the equations for $h_{ij}$, $i,j>1$. Although the equation that they satisfy seems 
to be 
the simplest of the ones that we have looked at, the $h_{ij}$'s are in fact the most subtle case. As we mentioned above, the issue is that the 
equation $Ah_{ij}=0$ is an Euler equation with no 0th order term. Therefore, there are constant solutions, which neither blow up nor go to zero as $r\rightarrow\infty$.

Since we want to invoke the rate at which $Lh_{i}$ tends toward 0, we will now be working with the $h_{i}$'s rather 
than $h$. Thus, we will need to quantify the rate at which the $h_{i}$'s are converging to their $T^{n-1}$-invariant 
limit.

Let 
\begin{equation}
\hat h_{i}(r) = \frac{1}{A(T(r))} \int_{T(r)} h_{i}(r, x) dA 
\end{equation}
be the average of $h$ over the torus $T(r)$. Then we have 
\begin{proposition}
$\|h_{i}- \hat h_{i}(r)\|_{C^{2}} = O\left(\frac{r}{R_{i}}\right)$
\end{proposition}
\textbf{Proof:} We have 
\begin{equation}
\sup_{x \in T(r)} |h_{i}(r,x) - \hat h_{i}(r)| \leq   \frac{1}{A(T(r))}\int_{T^{n-1}(r)} |h_{i}(r, x) - \hat h_{i}(r)| \ dA 
\end{equation}
Now, we know that $h_{i}$ is the lift of a form defined on a torus of diameter $O\left(\frac{r}{R_{i}}\right)$. Since 
$\|h_{i}\|_{k,\alpha}\leq C$, we know that the integrand must be less than $C$ times the diameter of the base torus, so 
it is also $O\left(\frac{r}{R_{i}}\right)$. We have assumed that $k\geq3$, so we can repeat this for the first and 
second derivatives of $h_{i}$. \qed

We know that on our unwrapped black hole metrics, we have 
\begin{equation}
L(\hat h) = L_{C}(\hat h) + O(r^{-(n-1)})
\end{equation}
 
Finally, since 
$\|Lh\|_{L^{\infty}_{r_{c}}}=\|\phi_{c}^{-1}Lh\|_{L^{\infty}}\rightarrow 0$,
\begin{equation}
L(h) =  o\left(\phi_{c}\right) 
\end{equation}

Putting this all together gives us 
\begin{eqnarray}
L_{C}(\hat h_{i}) & = & L_{BH}(\hat h_{i}) + O\left(r^{-(n-1)}\right) \\
  & = &   L_{BH}(h_{i}) + O\left(r^{-(n-1)}\right) + O\left(\frac{r}{R_{i}}\right)  \\
 & = &  o\left(\phi_{c}\right) + O\left(r^{-(n-1)}\right) + O\left(\frac{r}{R_{i}}\right)
\end{eqnarray}

Now, dropping the $i$'s and the hats, we see that if $a,b>1$
\begin{equation}
r^{2}h_{ab}^{\prime \prime} + nrh^{\prime}_{ab} = e_{ab}
\end{equation}
where 
\begin{equation}
e_{ab} =  O\left(\frac{r}{R_{i}}\right) + O\left(r^{-(n-1)}\right) +  o\left(\phi_{c}\right)
\end{equation}

Recall that we know that $\lim_{i\rightarrow \infty}h_{i}(r_{i})=0$. 
We want to show that this is true for any sequence 
$\rho_{i}\leq r_{i}$ with $\rho_{i}\rightarrow \infty$.

Using an integrating factor, we get  
\begin{eqnarray}
|h_{ab}(r_{i})-h_{ab}(\rho_{i})| & \leq & \int_{\rho_{i}}^{r_{i}}\frac{1}{r^{n}}\int_{r+}^{r} |e_{ab}(s)| s^{n-2}\ ds dr 
\\
&  \leq & \int_{\rho_{i}}^{r_{i}}\frac{1}{r^{n}}\left( C_{1}\int_{r+}^{r} \frac{s^{n-1}}{R_{i}}\ ds  +  
C_{2}\int_{r+}^{r} s^{-1}\ ds \right. \\
& & \ \ \left. + c_{i}  \int_{r+}^{r}s^{n-2}\phi_{c}(s)\  ds\right) dr
\end{eqnarray}
where $c_{i} \rightarrow 0$. Then this is 
\begin{eqnarray}
& \leq & C_{3}\frac{r_{i}}{R_{i}} + C_{4}r_{i}^{-(n-1)} + c_{i}\frac{r_{i}}{R_{i}}
\end{eqnarray}

Thus, we can conclude that 
\begin{equation}
\lim_{r \rightarrow \infty} |h(r)| = 0
\end{equation}

So by Corollary \ref{defs} we finally have our contradiction, and therefore the main 
estimate. \qed 

The reader may be somewhat uneasy at the following aspect of the above proof: the fact that the distance between a sequence of pairs of points could grow to infinity within a filling region caused us some difficulty at getting a uniform control on $h$ on the entire filling region. Should this problem not become even more serious when comparing points in different filling regions? The answer is no, since the weight functions are defined locally on each filling region, and are fine-tuned to the size of each one.

Finally, to finish this chapter, we extend the invertibility of $L_{g}= 2D_{g}\Phi_{g_{\sigma}}$ to a neighborhood of our approximate solution. Below, $B(x,\epsilon)$ will refer to a ball in the $S^{k, \alpha}_{\delta, r_{c}}$-topology, for the same choices of $k, \alpha, \delta$ and $r_{c}$ as above.

\begin{proposition}
There exist $\epsilon>0, \Lambda>0$ such that for all $\sigma$ large enough, we can choose $r_{c}$ such that the operator $L_{g}$ is 
invertible on the ball $B(g_{\sigma}, \epsilon)$, and for all $f \in \Phi(B(g_{\sigma}, \epsilon))$, we have that
\begin{equation}
\|(L_{g})^{-1}f\|_{k,\alpha, \delta, r_{c}}\leq \Lambda \|f\|_{k-2,\alpha,\delta,r_{c}}
\end{equation}
\end{proposition}
\textbf{Proof:} If not, there is a sequence of $g_{i}$'s and $\sigma_{i}$'s with $\|g_{i}-g_{\sigma_{i}}\|_{k, \alpha, 
\delta, r_{c}} \rightarrow 0$, and a sequence $h_{i} \in S^{k,\alpha}_{\delta, r_{c}}(M_{i},g_{i})$ such that
\begin{equation}
\|h_{i}\|_{k,\alpha,\delta, r_{c}}=1
\end{equation}
but 
\begin{equation}
\|L_{g_{i}}h_{i}\|_{k-2,\alpha,\delta, r_{c}} \longrightarrow 0
\end{equation}
But then we can repeat the proof of the previous proposition to obtain a contradiction. \qed

\section{Conclusions}

In this section, we conclude the construction of the AHE metrics on the $M_{\sigma}$'s.
\begin{proposition}
If $|\sigma|$ is large enough, then the manifold $M_{\sigma}$ admits an 
asymptotically hyperbolic Einstein manifold with the same conformal infinity as $N$.
\end{proposition}
\textbf{Proof:} Let $\delta \in \left(\frac{n-1}{2}, n-1\right)$. There is some $\epsilon>0$ such that for all $g$ in the ball 
$B(g_{\sigma}, \epsilon)$, the map
\begin{equation}
\Phi_{g_{\sigma}}: S^{k,\alpha}_{\delta, r_{c}}(g) \longrightarrow S^{k-2,\alpha}_{\delta, r_{c}}(g)
\end{equation}
has an invertible linearization and 
\begin{equation}
\|(D\Phi)^{-1}f\|_{k, \alpha, \delta,r_{c}} \leq \Lambda\|f\|_{k-2, \alpha,\delta, r_{c}}
\end{equation}
Thus, by the inverse function theorem(\cite{30}), $\Phi$ is invertible on $B(g_{\sigma},\epsilon)$, and maps $B(g_{\sigma},\epsilon)$ surjectively onto some $\mathcal{U}\subset S^{k-2,\alpha}_{\delta, r_{c}}$ containing $\Phi_{g_{\sigma}}(g_{\sigma})$. Thus, we need to show that $0 \in \mathcal{U}$. To do this, we will need a lower bound on the diameter of $\mathcal{U}$.

Let $B(\Phi(g_{\sigma}), \gamma) \subset \mathcal{U}$. By our control on $(D\Phi)^{-1}$ we know that $\Phi^{-1}$ is Lipschitz with Lipschitz constant $\Lambda$. Therefore,
\begin{equation}
\Phi^{-1}(B(\Phi(g_{\sigma}), \gamma)) \subseteq B(g_{\sigma}, \Lambda\gamma)
\end{equation}

Thus if we choose $\gamma=\frac{\epsilon}{\Lambda}$, we will obtain that 
$B(\Phi(g_{\sigma}),\gamma)\subseteq Im(\Phi)$. All that is left to do to show the existence of a solution to $\Phi=0$ is to make sure that $0\in 
B(\Phi(g_{\sigma}),\gamma)$ for $\sigma$ large enough. But 
$\|\Phi(g_{\sigma})\|_{k-2,\alpha,\delta,r_{c}} = 
O(|\sigma|^{1-n})$, so \begin{equation}
\|\Phi(g_{\sigma}) - 0\|_{k-2,\alpha, \delta,r_{c}}\leq C|\sigma|^{1-n} \leq \frac{\epsilon}{\Lambda} 
\end{equation}
for $|\sigma|$ large enough, since $\epsilon$ and $\Lambda$ are fixed.

Let $g_{E}=\Phi^{-1}(0)$. Since $g_{\sigma}$ has negative Ricci curvature, we may assume that $\epsilon$ is small enough that $g_{E}$ does too. Since $B_{g_{\sigma}}(g_{\sigma})= 0$, and $\lim_{r\rightarrow\infty}\|g-g_{\sigma}\|_{k,\alpha, \delta, r_{c}}=0$ we can conclude that $\lim_{r \rightarrow \infty} \left|B_{g_{\sigma}}(g_{E})\right| = 0$. Invoking lemma \ref{Biq} gives us that $g_{E}$ is an Einstein metric. Since $g_{E}$ is a perturbation of $g_{\sigma}$ by an element of $S^{k, \alpha}_{\delta, r_{c}}$, it must have the same conformal infinity, so we are done. \qed

Note that by construction, all of these metrics are nondegenerate, i.e. $D\Phi$ has no $L^2$ kernel. Furthermore, all of these metrics are isolated points in 
the moduli space of AH Einstein metrics on $M_{\sigma}$ with boundary metric $[\gamma]$, equipped with the $C^{k,\alpha}_{\delta}$ topology.               

Abusing notation slightly, let us denote these AHE metrics by $g_{\sigma}$. Then we have that, for any sequence of points $p_{\sigma}$ which remain within a bounded distance of a gluing torus
\begin{equation}
\lim_{|\sigma|\rightarrow \infty}(M_{\sigma},g_{\sigma},p_{\sigma}) = (N,g)
\end{equation}
in any of the $S^{k,\alpha}_{\delta}$ topologies, where $(N,g)$ is our original hyperbolic manifold. By only opening one cusp at a time, we will construct nonhyperbolic AHE manifolds with cusps. 

\textbf{Proof of Theorem \ref{open}:}
Let $\sigma = (\sigma_{1}, ..., \sigma_{l})$ have large enough norm that the
manifold $M_{\sigma}$ admits an Einstein metric. Then define a sequence $\sigma_{i} =
(\upsilon_{i}, \sigma_{2}, ..., \sigma_{l})$, where $\upsilon_{i}$ is a sequence of
geodesics such that $L(\upsilon_{i}) \rightarrow \infty$. Let $p_{i}$ be a sequence of
points within a bounded distance of a given gluing torus. Then
\begin{equation}
(M,g_{\infty}) = \lim_{i\rightarrow \infty}(M_{i},g_{i},p_{i})
\end{equation}
is a nonhyperbolic AHE manifolds with cusps. The conformal infinity of
$(M,g_{\infty})$ is the same as that of $(N,g)$. \qed.

All of the conformal infinities involved in our constructions are 
necessarily conformally flat, since they are the conformal infinity 
of a hyperbolic manifold. All of the examples known to the author of 
conformal classes bounding infinitely many AHE manifolds are conformally flat.

It is possible to extend this filling construction to some other types of finite-volume cusp ends. In such a case, the 
cusp cross-section is necessarily a compact flat manifold, and thus a 
finite quotient of a torus by Bieberbach's 
theorem. One can find sufficient conditions to construct an approximate 
solution in terms of the geometry of the flat cross-section. In all 
dimensions greater than 3 there exist flat cusp cross-sections to which 
this procedure cannot be 
applied. For more on this, see (\cite{1}).

\end{document}